\documentclass[11pt,reqno]{amsart}
\usepackage{amssymb}
\usepackage{amsfonts}
\usepackage{amsthm}
\usepackage{amscd}
\usepackage{mathptmx}       
\usepackage{helvet}         
\usepackage{courier}        
\usepackage{type1cm}        
\usepackage{graphicx}        
\usepackage{amsmath}
\usepackage{amssymb}
\usepackage{amsfonts}
\usepackage{amsthm}
\usepackage{amscd}
\theoremstyle{plain}
\newtheorem{theorem}{Theorem}[section]
\newtheorem{proposition}[theorem]{Proposition}

\theoremstyle{definition}
\newtheorem{definition}[theorem]{Definition}
\newtheorem{remark}[theorem]{Remark}
\newtheorem*{example}{Example}
\newcommand{\refE}[1]{(\ref{eq:#1})}

\newcommand{\refS}[1]{Section~\ref{sec:#1}}

\newcommand{\refD}[1]{Definition~\ref{D:#1}}

\newcommand{\C}{\ensuremath{\mathbb{C}}}

\newcommand{\Z}{\ensuremath{\mathbb{Z}}}
\newcommand{\K}{\ensuremath{\mathcal{K}}}

\newcommand{\J}{\ensuremath{\mathbb{J}}}
\newcommand{\Jl}{\ensuremath{\mathbb{J}_\lambda}}
\newcommand{\Pro}{\ensuremath{\mathbb{P}}}

\newcommand{\tr}{\mathrm{tr}}

\newcommand{\cins}{\frac 1{2\pi\mathrm{i}}\int_{C_S}}

\newcommand{\cint}[1]{\frac 1{2\pi\mathrm{i}}\int_{#1}}

\newcommand{\cintl}[1]{\frac 1{24\pi\mathrm{i}}\int_{#1 }}
\newcommand{\g}{\mathfrak{g}}
\newcommand{\gb}{\overline{\mathfrak{g}}}
\newcommand{\gh}{\widehat{\mathfrak{g}}}

\newcommand{\A}{\mathcal{A}}
\newcommand{\W}{\mathcal{W}}
\newcommand{\V}{\mathcal{V}}
\newcommand{\U}{\mathcal{U}}
\newcommand{\Ja}{\mathcal{J}}

\newcommand{\La}{\mathcal{L}}
\newcommand{\Sp}{\mathcal{F}^{-1/2}}
\newcommand{\Sa}{\mathcal{S}}

\newcommand{\Lh}{\widehat{\mathcal{L}}}

\newcommand{\ord}{\operatorname{ord}}
\newcommand{\res}{\operatorname{res}}
\newcommand{\fpz}{\frac {d }{dz}}
\newcommand{\pfz}[1]{\frac {d#1}{dz}}

\newcommand{\Ho}{\mathrm{H}}

\newcommand{\al}{\ensuremath{\alpha}}
\newcommand{\be}{\ensuremath{\beta}}

\newcommand{\Hl}[1][\lambda]{\mathcal{H}^{#1}}
\newcommand{\Fl}[1][\lambda]{\mathcal{F}^{#1}}

\newcommand{\Dal}{\mathcal{D}_\lambda}
\newcommand{\Dalh}{\widehat{\mathcal{D}}_\lambda}
\newcommand{\Do}{\mathcal{D}^1}
\newcommand{\Dh}{\widehat{\mathcal{D}^1}}
\newcommand{\kndual}[2]{\langle #1,#2\rangle}
\newcommand{\ldot}{\,.\,}

\newcommand{\glih}{\widehat{gl}(\infty)}
\newcommand{\glib}{\overline{gl}(\infty)}
\newcommand{\la}{\lambda}
\newcommand{\ka}{\kappa}

\newcommand{\sbul}{{\boldsymbol{\cdot}}}

\newcommand{\ddx}[1][x]{\frac{d}{d #1}}
\newcommand{\ddz}[1][z]{\frac{d}{d #1}}
\newcommand{\zinf}{\{0,\infty\}}
\newcommand{\set}[1]{\left\{#1\right\}} 

\newcommand{\sln}{\mathfrak{sl}}

\newcommand{\so}{\mathfrak{so}}
\newcommand{\spn}{\mathfrak{sp}}
\newcommand{\gl}{\mathfrak{gl}}

\newcommand{\ga}{{\gamma}}

\begin{document}
\title[Krichever-Novikov type algebras]
{From the Virasoro Algebra to 
\\ Krichever--Novikov Type 
Algebras 
\\
and Beyond}
\author[Martin Schlichenmaier]{Martin Schlichenmaier}
\thanks{Partial  support by
the 
Internal Research Project  GEOMQ11,  University of Luxembourg,
is acknowledged.}
\address{%
University of Luxembourg\\
Mathematics Research Unit, FSTC\\
Campus Kirchberg\\ 6, rue Coudenhove-Kalergi,
L-1359 Luxembourg-Kirchberg\\ Luxembourg
}
\email{martin.schlichenmaier@uni.lu}
\begin{abstract}
Starting from the Virasoro algebra and its relatives the
generalization
to higher genus compact Riemann surfaces was initiated by Krichever and 
Novikov. The elements of these algebras are meromorphic objects 
which are holomorphic outside a finite set of points.
A crucial and non-trivial point is to establish an almost-grading
replacing the honest grading in the Virasoro case. 
Such an almost-grading is given by splitting the set of points of 
possible  poles into two non-empty disjoint subsets.
Krichever and Novikov considered the two-point case. Schlichenmaier
studied the most general multi-point situation with arbitrary splittings.
Here we will review the path of developments from the Virasoro
algebra to  its higher genus and multi-point analogs.
The starting point will be a Poisson algebra 
structure on the space of meromorphic forms of all weights. 
As sub-structures the vector field algebras, function algebras, 
Lie superalgebras and the related current algebras show up.
All these algebras will be almost-graded.
In detail almost-graded central extensions are classified.
In particular, for the vector field algebra it is
essentially unique. 
The defining cocycle are given in geometric terms.
Some applications, including the semi-infinite wedge form
representations
are recalled. Finally, some remarks on the 
by Krichever and Sheinman recently introduced 
Lax operator algebras are made.
\end{abstract}

\maketitle
\section{Introduction}\label{sec:intro}

Lie groups and Lie algebras are related to 
symmetries of systems.
By the use of the symmetry the system can be better understood, 
maybe  it is even possible to solve it in a certain sense.
Here we deal with systems which have an infinite number of
independent degrees of freedom. They appear for example in 
Conformal Field Theory (CFT),  see e.g. \cite{BP}, \cite{TUY}.
 But also in the theory of 
partial differential equations and at 
many other places in- and outside of
mathematics they play an important role. 
The appearing Lie groups and Lie algebras 
are infinite dimensional.
Some of the  simplest  nontrivial 
infinite dimensional Lie algebras are the  Witt algebra 
and its central extension the Virasoro algebra.
We will recall their definitions in \refS{vir}.
In the  sense explained (in particular in CFT)
they are related to what is called the genus zero situation.
For CFT on arbitrary genus Riemann surfaces the
Krichever-Novikov (KN) type algebras, to be discussed here,
will show up as algebras of 
global symmetry operators.

These algebras are defined 
via meromorphic objects on compact Riemann surfaces $\Sigma$ of
arbitrary genus with controlled polar behaviour. More precisely,
poles are only allowed at a fixed finite set of points denoted
by $A$. The ``classical'' examples are the algebras defined by objects
on the Riemann sphere (genus zero) with possible poles only at
$\{0,\infty\}$. This yields e.g. the  Witt algebra, 
the classical current
algebras,
including  their central extensions the Virasoro, and the affine Kac-Moody
algebras \cite{KacB}.
For higher genus, but still only for two points where poles are
allowed, they were generalised by Krichever and Novikov 
\cite{KNFa}, \cite{KNFb}, \cite{KNFc} in 1987.
In 1990 the author 
\cite{Schlmp}, \cite{Schleg}, \cite{Schlce}, \cite{SchlDiss} extended the approach
further to the general multi-point case. 

This extension was not a straight-forward generalization.
The crucial point is to introduce a replacement of the
graded algebra structure present in the ``classical'' case.
Krichever and Novikov found that an almost-grading, see
\refD{almgrad} below, will be enough  
to do the usual constructions in
representation theory, like triangular decompositions, highest
weight modules, Verma modules which  are demanded by the applications.
In \cite{Schlce}, \cite{SchlDiss} it was realized that a splitting
of $A$ into two disjoint non-empty subsets $A=I\cup O$ is crucial for
introducing an almost-grading
and the corresponding almost-grading was given.
In the two-point situation  there is only one such splitting
(up to inversion) hence there is only one almost-grading, which 
in the classical case is a honest grading.
Similar to  the classical situation a Krichever-Novikov algebra,
should
always be considered as an algebra of meromorphic objects
with an almost-grading coming from such a fixed splitting.

I like to point out that already in the genus zero case 
(i.e. the Riemann sphere case) with more
than two points where  poles are allowed the algebras will only be 
almost-graded. In fact, quite a number of interesting
new phenomena will show up already there, see
\cite{SchlDeg}, \cite{FiaSchl1}, \cite{FiaSchlaff}, \cite{Brem2}.

In this review no proofs  are supplied.
For them I have to refer to 
the original articles and/or to the forthcoming book
\cite{Schlknbook}. For some applications jointly obtained
with Oleg Sheinman, see also 
\cite{SheBo}.
For more on the Witt and Virasoro algebra see for example
the book \cite{GR}.

\medskip

After recalling the definition of the Witt and Virasoro algebra in
\refS{vir} we start with describing the geometric set-up of 
Krichever-Novikov (KN) type algebras in \refS{algebra}.
We introduce a Poisson algebra structure on the space of
meromorphic forms (holomorphic outside of the fixed  set $A$ of 
points where poles are allowed) of all weights  (integer and
half-integer).
Special substructures will yield the function algebra, the vector
field
algebra and more generally the differential operator
algebra. Moreover,
we discuss also the Lie superalgebras of KN type defined via
forms of weight -1/2. An important example role also is played
by the current algebra (arbitrary genus - multi-point) 
associated to a finite-dimensional Lie algebra. 

In \refS{almgrad} we introduce the almost-grading induced by the
splitting of $A$ into ``incoming'' and ``outgoing'' points,
$A=I\cup O$.

In \refS{centvec} we discuss central extensions for our algebras.
Central extensions appear naturally in the context of
quantization and regularization of actions. 
We give for all our algebras geometrically defined central extensions.
The defining cocycle for the Virasoro algebra obviously does not
make any sense in the higher genus and/or multi-point case. 
For the geometric description  we use 
projective and affine connections.
In contrast to the classical case there are a many inequivalent
cocycles
and central extensions. 
If we restrict our attention to the cases where we can extend the
almost-grading to the central extensions the author obtained complete
classification and uniqueness results. They are
described in \refS{class}.

In \refS{fr} we present further results. In particular,
we discuss how from the representation of the vector field
algebra (or more general of the differential operator algebra) on the
forms
of weight $\lambda$  one obtains semi-infinite wedge representations
(fermionic Fock space representations) of the centrally extended 
algebras.
These representations have ground states (vacua), creation 
and annihilation operators.
We add some words about $b-c$ systems, Sugawara construction,
Wess-Zumino-Novikov-Witten (WZNW) models, Knizhnik-Zamolodchikov
(KZ) connections, and deformations of the Virasoro algebra.

Recently, a new class of current type algebras 
the
Lax operator algebras,
were introduced by Krichever and Sheinman
\cite{Klax}, \cite{KSlax}.
I will report on them in \refS{lax}.

In the closing \refS{hist}
 some historical remarks (also on related works)
on 
Krichever-Novikov type algebras
and some references are given. More references can be found in
\cite{Schlknbook}.


\section{The Witt and  Virasoro Algebra}
\label{sec:vir}
\subsection{The Witt Algebra}

The {\it Witt algebra} $\W$, also sometimes called 
Virasoro algebra without central
  term%
\footnote{
In the book  \cite{GR} arguments are given why it is more appropriate
just to use Virasoro algebra, as Witt introduced ``his'' algebra in
a characteristic $p$ context. 
Nevertheless, I decided to stick here to the most common
convention.},
is the complex Lie algebra  
generated as vector space by the
elements 
\ 
$\{e_n\mid n\in \mathbb{Z}\}$\  with 
Lie structure  
\begin{equation}\label{eq:Wstruct}
[e_n,e_m]=(m-n)e_{n+m},\quad n,m\in\Z\;.
\end{equation}
One of its realization is 
as complexification of the Lie algebra of
polynomial vector fields $Vect_{pol}(S^1)$ on the circle $S^1$,
which is a subalgebra of $Vect(S^1)$, 
die Lie algebra of all $C^\infty$ vector
fields  on the circle. 
In this realization 
\begin{equation} 
e_{n}:= -\mathrm{i}\,
\exp{\mathrm{i}\, n\varphi}\,\ddx [\varphi], \qquad  n\in\Z \;.
\end{equation}
The Lie product is the usual Lie bracket of vector fields.

If we extend these generators to 
 the whole punctured complex plane we obtain 
\begin{equation}
e_n=z^{n+1}\frac {d}{dz}, \qquad  n\in\Z \;.
\end{equation} 
This gives another realization of the Witt algebra 
as the 
algebra of those meromorphic vector fields on the Riemann sphere
$\Pro^1(\C)$ which are holomorphic outside $\{0\}$ and 
$\{\infty\}$. 

Let $z$ be the (quasi) global coordinate $z$ 
(quasi, because it is not defined at
$\infty$). 
Let $w=1/z$ be the local coordinate at $\infty$.
A global meromorphic vector field $v$ on  $\Pro^1(\C)$ will be given 
on the corresponding subsets where $z$ resp. $w$ are defined
as
\begin{equation}
v=\left(v_1(z)\ddz\ ,\  v_2(w)\ddx[w]\right),
\qquad v_2(w)=-v_1(z(w))w^2.
\end{equation}
The function $v_1$ will determine the 
vector field $v$. Hence, we will usually just 
write  $v_1$ and in fact identify the vector field
$v$ with its local representing function $v_1$, which we will
denote by the same letter.

For the bracket we calculate
\begin{equation}
[v,u]=\left(v\ddz u-u\ddz v\right)\ddz.
\end{equation}
The space of all meromorphic vector fields constitute a Lie algebra.
The subspace of those meromorphic vector fields which are holomorphic
outside of  $\set{0,\infty}$ is a Lie subalgebra.
Its  elements can be given as 
\begin{equation}
v(z)=f(z)\ddz
\end{equation}
where $f$ is a meromorphic function on $\Pro^1(\C)$, which is holomorphic
outside $\zinf$. Those are exactly the 
Laurent polynomials
$\C[z,z^{-1}]$.
Consequently, this subalgebra has  the set
$\ \{e_n, n\in\Z \}\ $ as basis 
elements. The Lie product is the same and it can be 
identified with the Witt algebra $\W$.

The subalgebra of global holomorphic vector fields
is ${\langle e_{-1},e_0,e_1\rangle}_\C$. It is
isomorphic to the Lie algebra $\mathfrak{sl}(2,\C)$.

The algebra 
$\W$ is more than just a Lie algebra. It is a graded Lie algebra.
If we set for the degree 
$\deg(e_n):=n$ 
then $\deg([e_n,e_m])=\deg (e_n)+\deg(e_m)$
and we obtain the degree decomposition
\begin{equation}
\W=\bigoplus_{n\in\Z}\W_n, \qquad 
\W_n={\langle e_n\rangle}_\C\; .
\end{equation}
Note that $[e_0,e_n]=n\,e_n$, which says that
the degree decomposition is the eigen-space decomposition with
respect to the adjoint action of $e_0$ on $\W$.

Algebraically $\W$ can also be given as
Lie algebra of derivations of the algebra of Laurent
polynomials $\C[z,z^{-1}]$.

\subsection{The Virasoro Algebra}

In the process of quantizing or regularization  
one is often forced to
modify an action of a Lie algebra. 
A typical example is given by the product of infinite sums of
operators. Quite often they are only well-defined if a certain
``normal ordering'' is introduced. 
In this way  the modified action will only be 
a projective action. This can be made to an honest Lie action by
passing  to a suitable central extension of the Lie algebra.

For the Witt algebra the  universal one-dimensional
central extension is the 
{\it Virasoro algebra} $\V$. 
As vector space it is the direct sum $\V=\C\oplus \W$.
If we set for $x\in\W$, 
$\hat x:=(0,x)$, and $t:=(1,0)$
then 
its basis elements are $\hat e_n, \ n\in\Z$ and 
$t$ with the  Lie product
\begin{equation}\label{eq:Vstruct}
[\hat e_n,\hat e_m]=(m-n)\hat e_{n+m}-\frac 1{12}(n^3-n)\delta_n^{-m}\,t,\quad 
\quad
[\hat e_n,t]=[t,t]=0\;,
\end{equation}
for%
\footnote{
Here $\delta_k^l$ is the Kronecker delta which is equal to 1 if
$k=l$, otherwise zero.}
all 
 $n,m\in\Z$.
If we set $\deg(\hat e_n):=\deg(e_n)=n$ and $\deg(t):=0$ then 
$\V$ becomes a graded algebra. 
The algebra $\W$ will only be a subspace, not a subalgebra 
of $\V$. It will be a quotient.
In some abuse of notation we identify the element $\hat x\in\V$ with
$x\in\W$.
Up to equivalence and rescaling the central element $t$, this is beside
the trivial (splitting) central extension the only central extension.

\section{The Krichever-Novikov Type Algebras}\label{sec:algebra}
\subsection{The Geometric Set-Up}
For the whole article let $\Sigma$ be a compact Riemann surface 
without any restriction for the genus $g=g(\Sigma)$.
Furthermore, let $A$ be a finite subset of $\Sigma$.
Later we will need a splitting of $A$ into two non-empty disjoint
subsets $I$ and $O$, i.e. $A=I\cup O$. Set $N:=\#A$,
$K:=\#I$, $M:=\#O$, with $N=K+M$. 
More precisely, let
\begin{equation}
I=(P_1,\ldots,P_K),\quad\text{and}\quad
O=(Q_1,\ldots,Q_{M})
\end{equation}
be disjoint  ordered tuples of  distinct points (``marked points'',
``punctures'') on the Riemann surface.
In particular, we assume $P_i\ne Q_j$ for every
pair $(i,j)$. The points in $I$ are
called the {\it in-points}, the points in $O$ the {\it out-points}.
Sometimes we consider $I$ and $O$ simply as sets.

In the article we sometimes refer to the classical situation. By this
we understand 
\begin{equation}\label{eq:class}
\Sigma=\Pro^1(\C)=S^2, \quad I=\{z=0\},\quad 
 O=\{z=\infty\}
\end{equation}

The following Figures 
\ref{f3}, \ref{f2}, \ref{f1} 
exemplify the different situations:
%
\begin{figure}[h]
\centering{\includegraphics[width=0.5\textwidth]
{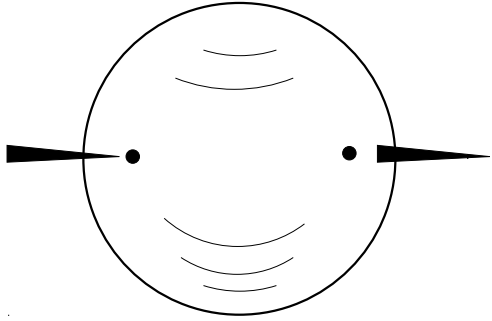}}
\vspace*{8pt}
\caption{Riemann surface of genus zero 
 with one incoming and one outgoing point. 
\label{f3}}
\end{figure}
%
%
\begin{figure}[h]
\centering{\includegraphics[width=0.5\textwidth]
{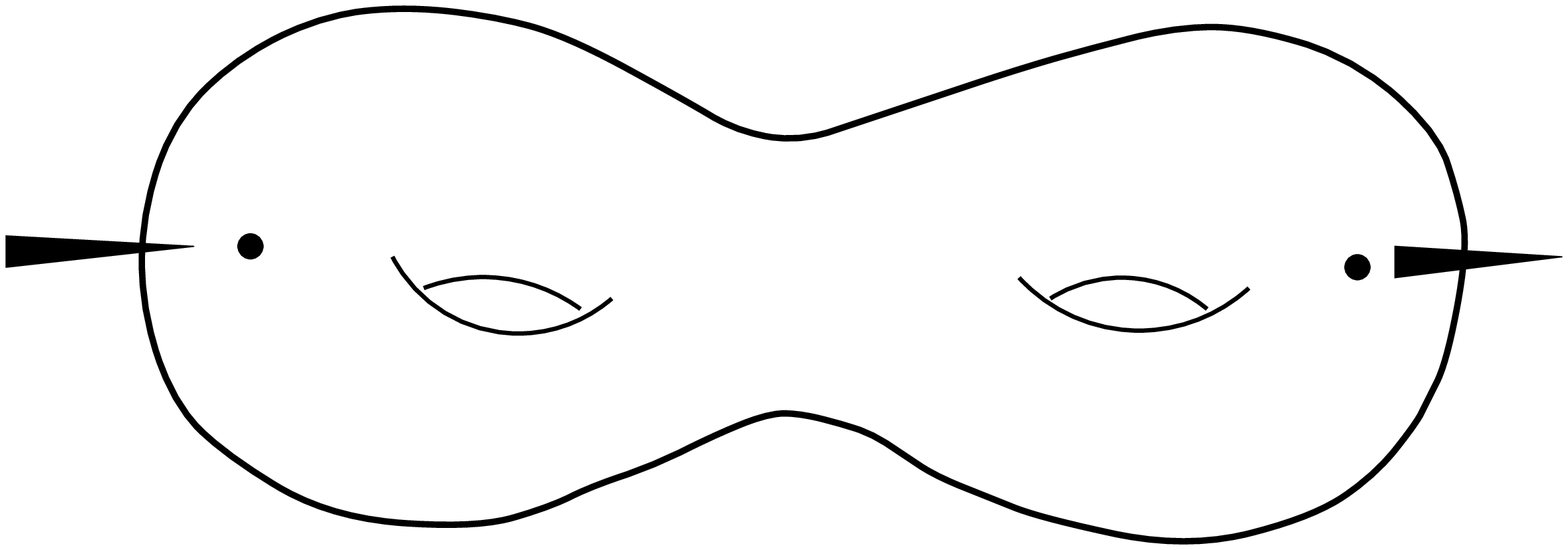}}
\vspace*{8pt}
\caption{Riemann surface of genus two 
 with one incoming and one outgoing point. 
\label{f2}}
\end{figure}
%
\begin{figure}[h]
\centerline{\includegraphics[width=0.5\textwidth]
{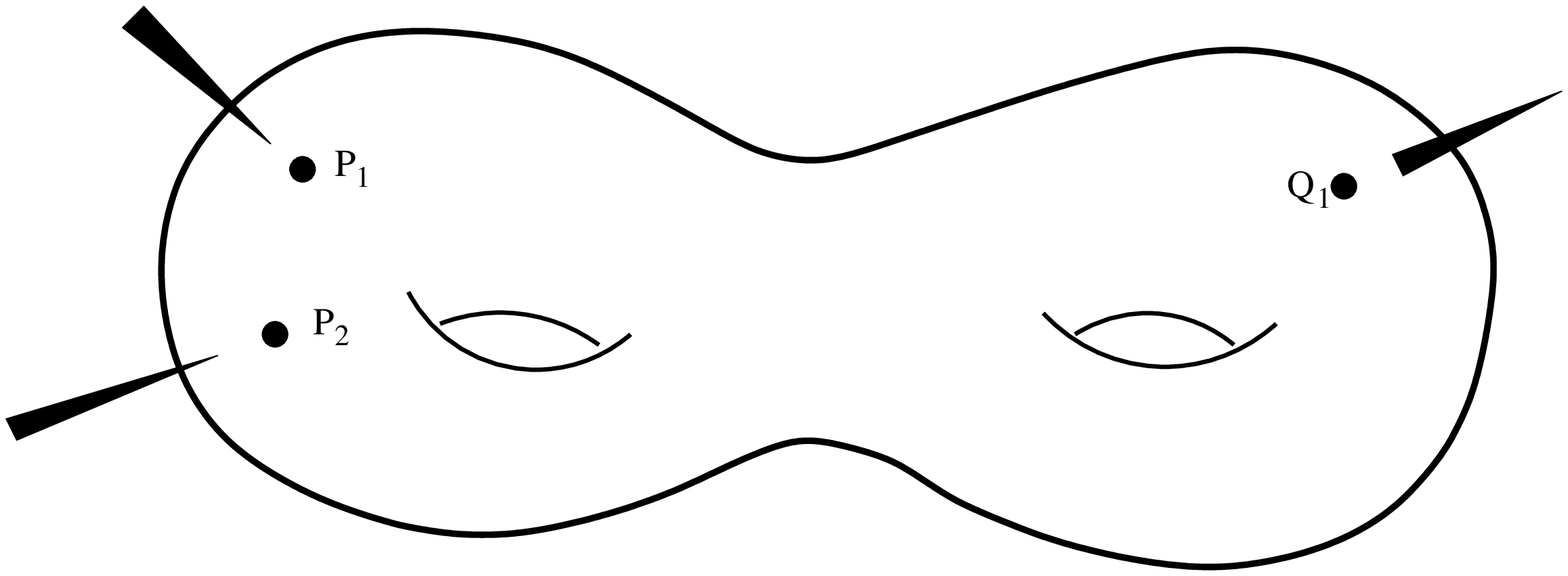}}
\vspace*{8pt}
\caption{Riemann surface of genus two 
 with two incoming points and one outgoing point. 
\label{f1}}
\end{figure}

Our objects, algebras, structures, ...  will be
meromorphic objects 
defined on $\Sigma$ which are holomorphic outside of the points in
$A$.
To introduce the objects 
let $\K=\K_\Sigma$ be the canonical line bundle of $\Sigma$,
resp. the locally free canonically sheaf.
The local sections of the bundle are
the local holomorphic differentials.
If $P\in\Sigma$ is a point and $z$ a local holomorphic coordinate
at $P$ then a local holomorphic differential 
can be written as $f(z)dz$ with a local holomorphic function 
$f$ defined in a neighbourhood of $P$.
A global holomorphic section can be described locally with respect to
 a covering
by coordinate charts $(U_i,z_i)_{i\in J}$ by a system of
local holomorphic functions $(f_i)_{i\in J}$, which are related by 
the transformation rule induced by the
coordinate change map $z_j=z_j(z_i)$ and the condition
$f_idz_i=f_jdz_j$ yielding
\begin{equation}\label{eq:ktrans}
f_j=f_i\cdot \left(\frac {dz_j}{dz_i}\right)^{-1}.
\end{equation}
Moreover, 
a meromorphic section of $\K$ is given as 
a collection of local meromorphic functions $(h_i)_{i\in J}$ 
for which the transformation
law \refE{ktrans} still is true.

In the following $\la$ is either an integer or a half-integer.
If $\la$ is an integer then 
\newline
(1) $\K^{\la}=\K^{\otimes \la}$ 
for $\la>0$, 
\newline
(2)  $\K^{0}=\mathcal{O}$,
 the trivial line bundle, and 
\newline
(3) $\K^{\la}=(\K^*)^{\otimes (-\la)}$ 
for $\la<0$.
\newline
Here as usual $\K^*$ denotes the  dual line bundle
to the canonical line bundle.
The dual line bundle is the holomorphic tangent line bundle, whose
local sections are the holomorphic tangent vector fields
$f(z)(d/d z)$.
If $\la$ is a half-integer, then we first have to fix a ``square
root''
of the canonical line bundle, sometimes called a
\emph{theta-characteristics}.
This means we fix a line bundle $L$ for
which $L^{\otimes 2}=\K$.

After such a choice of $L$ is done we set 
$\K^{\la}=\K^{\la}_L=L^{\otimes 2\la}$.
In most cases we will drop the mentioning of $L$, but we have to
keep the choice in mind.
Also the fine-structure of the 
algebras we are about to define will depend on the
choice. But the main  properties will remain the same.

\begin{remark}
A Riemann surface of genus $g$ has exactly 
$2^{2g}$ non-isomorphic square roots of $\K$. 
For $g=0$ we have $\K=\mathcal{O}(-2)$, and 
$L=\mathcal{O}(-1)$, the tautological bundle,
 is the unique square root.
Already for $g=1$ we have 4 non-isomorphic ones.
As in this case $\K=\mathcal{O}$ one solution is 
$L_0=\mathcal{O}$. 
But we have also other bundles $L_i$, $i=1,2,3$.
Note that $L_0$ has a non-vanishing global holomorphic
section, whereas this is not the case for 
$L_1,L_2$, $L_3$.
In general, 
depending on the parity of $\dim\Ho(\Sigma,L)$, one 
distinguishes even and odd theta characteristics $L$. For $g=1$ 
the bundle $\mathcal{O}$ is an odd, the others are even
theta characteristics.
\end{remark}

\bigskip
We set
\begin{multline}\qquad
\Fl:=\Fl(A):=\{f \text{ is a global meromorphic section of } K^\la
\mid
\\
\text{such that} \ 
f \text{ is  holomorphic over } \Sigma\setminus A\}\;.\qquad
\end{multline}
We will drop the set $A$ in
the notation.
Obviously, $\Fl$ is an infinite dimensional  $\C$-vector space.
Recall that in the case of half-integer $\la$ everything depends on
the theta characteristic $L$.

The elements of the space $\Fl$ we  call \emph{meromorphic forms
of weight $\la$} 
(with respect to the 
theta characteristic $L$).
In local coordinates $z_i$ we can write such a form as 
$f_idz_i^\la$, with $f_i$ a local holomorphic, resp. meromorphic form.

Special important cases of the weights are
the functions ($\la=0$), the space is also denoted by $\A$,
the vector fields ($\la=-1$), denoted by $\La$,
the differentials ($\la=1$), and the quadratic
differentials ($\la=2$).

\medskip
Next we introduce algebraic operations on the space of all
weights
\begin{equation}
\mathcal{F}:=\bigoplus_{\la\in\frac 12\Z}\Fl.
\end{equation}
These operations will allow us to introduce the 
algebras we are heading for.

\subsection{Associative Structure}
The natural map of the 
locally free sheaves of rang one  
\begin{equation}
\K^{\lambda}\times \K^\nu\to
 \K^{\lambda}\otimes \K^\nu\cong 
\K^{\lambda+\nu},
\quad
(s,t)\mapsto s\otimes t,
\end{equation}
defines a bilinear map 
\begin{equation}\label{eq:afl}
\sbul:\Fl\times \Fl[\nu]\to \Fl[\la+\nu].
\end{equation}
With respect to local trivialisations this corresponds to the
multiplication of the local representing meromorphic functions
\begin{equation}
(s\, dz^{\la},t\, dz^{\nu} )
\mapsto s\, dz^{\la}\;\sbul\;t\, dz^{\nu}= s\cdot t\;
dz^{\la+\nu}.
\end{equation}
If there is no danger of confusion then we will mostly
use the same symbol for the section and for the local
representing function.

The following is obvious
\begin{proposition}\label{P:aass}
The vector space 
$\mathcal{F}$ is an  associative and commutative
graded (over $\frac 12\Z$)
algebras. Moreover, $\A=\Fl[0]$ is a subalgebra.
\end{proposition}
\begin{definition}
The associative algebra $\A$ is the Krichever-Novikov function
algebra (associated to $(\Sigma,A)$).
\end{definition}
Of course, it is the  algebra of meromorphic functions on $\Sigma$
which are holomorphic outside of $A$.
The spaces $\Fl$ are modules over $\A$.
In the classical situation $\A=\C[z,z^{-1}]$, the
algebra of Laurent polynomials.
\subsection{Lie Algebra Structure}
Next we define a Lie  algebra structure on the space $\mathcal{F}$.
The structure is  induced by the map 
\begin{equation}\label{eq:lp}
\Fl\times \Fl[\nu]\to \Fl[\la+\nu+1],
\qquad (s,t)\mapsto [s,t],
\end{equation}
which is defined in local representatives of the sections by
\begin{equation}\label{eq:aliea}
(s\, dz^{\la},t\, dz^{\nu} )
\mapsto [s\, dz^{\la},t\, dz^{\nu}]:= \left((-\la)s\pfz{t}+\nu\, 
t\pfz{s}\right)
dz^{\la+\nu+1},
\end{equation}
and bilinearly extended to $\mathcal{F}$.
\begin{proposition} \cite{SchlHab}, \cite{Schlknbook}

\noindent
(a)
The bilinear map $[.,.]$ defines a Lie algebra structure 
on $\mathcal{F}$.

\noindent
(b) The space  $\mathcal{F}$
with respect to $\sbul$ and 
$[.,.]$ is a Poisson algebra.
\end{proposition}
Next we consider certain important substructures.
\subsection{The Vector Field Algebra and the Lie Derivative}
For  $\la=\nu=-1$ in \refE{lp} we end up in 
$\Fl[-1]$ again.
Hence, 
\begin{proposition}
The subspace $\La=\Fl[-1]$ is a Lie subalgebra, and the 
$\Fl$'s are Lie modules over $\La$.
\end{proposition} 
As forms of weight $-1$ are vector fields, $\La$ could also be
defined as the Lie algebra of those meromorphic vector fields 
on the Riemann surface $\Sigma$ which are holomorphic outside
of $A$.
The product \refE{aliea} gives the usual
Lie bracket of vector fields and the Lie derivative for their actions on 
forms. 
Due to its importance let us specialize this.
We obtain
(naming the  local functions
with the same symbol as the section) 
\begin{equation}\label{eq:aLbrack}
[e,f]_|(z)=[e(z)\fpz, f(z)\fpz]=
\left( e(z)\pfz f(z)- f(z)\pfz e(z)\right)\fpz \ ,
\end{equation}
\begin{equation}\label{eq:alied}
\nabla_e(g)_|(z)=L_e(g)_|=
e\ldot g_{|}=
\left( e(z)\pfz g(z)+\lambda g(z)\pfz e(z)\right)(dz)^{\la}\; .
\end{equation}
\begin{definition}
The algebra $\La$ is called \emph{Krichever-Novikov type vector
field algebra} (associated to $(\Sigma,A)$.
\end{definition}
In the classical case this gives the Witt algebra.

\subsection{The Algebra of Differential Operators}

In $\mathcal{F}$,
considered as Lie algebra,
$\A=\Fl[0]$ is an abelian Lie subalgebra and the vector space sum
$\Fl[0]\oplus\Fl[-1]=\A\oplus \La$  
is also a Lie subalgebra of  $\mathcal{F}$.
In an equivalent way it can also be constructed as 
semi-direct sum of $\A$ considered as abelian Lie algebra
and $\La$ operating on $\A$ by taking the derivative.
\begin{definition}
This Lie algebra is called the 
\emph{Lie algebra of differential operators of degree
$\le 1$} of KN type (associated to $(\Sigma,A)$)
and is denoted by
$\Do$.
\end{definition}
In more direct terms $\Do=\A\oplus \La$ as vector space direct sum
and endowed with the Lie product
\begin{equation}
[(g,e),(h,f)]=(e\ldot h-f\ldot g\,,\,[e,f]).
\end{equation} 
The spaces  $\Fl$ will be Lie-modules over $\Do$.

Its universal enveloping algebra will be the 
algebra of all differential operators of arbitrary degree
\cite{SchlDiss}, \cite{SchlDiff},
\cite{Schlloc}.


\subsection{The Superalgebra of Half Forms}
\label{sec:super}

Next we consider the associative product
\begin{equation}
\sbul\ \Fl[-1/2]\times \Fl[-1/2]\to \Fl[-1]=\La.
\end{equation}

We introduce the vector space and the product
\begin{equation}\label{eq:saprod}
\Sa:=\La\oplus\Fl[-1/2],\quad
[(e,\varphi),(f,\psi)]:=([e,f]+\varphi\;\sbul\;\psi,
e\ldot \varphi-f\ldot \psi).
\end{equation}
Usually we will denote
the elements of $\La$ by $e,f, \dots$, and the elements of
$\Sp$ by $\varphi,\psi,\ldots$.

The definition \refE{saprod} can be reformulated as an extension
of $[.,.]$ on  $\La$ to a ``super-bracket'' (denoted by the same
symbol) on $\Sa$ by setting
\begin{equation}\label{eq:sadef}
[e,\varphi]:=-[\varphi,e]:=e\ldot \varphi
=(e\frac {d\varphi}{dz}-\frac 12\varphi\frac {de}{dz})(dz)^{-1/2}
\end{equation}
and
\begin{equation}
[\varphi,\psi]:=\varphi\;\sbul\; \psi\; .
\end{equation}
We call the elements of $\La$ elements of even parity,
and the elements of $\Sp$ elements of odd parity.
For such elements $x$ we denote by $\bar x\in\{\bar 0,\bar 1\}$ their
parity. 

The sum \refE{saprod} can also  be described
as $\Sa=\Sa_{\bar 0}\oplus\Sa_{\bar 1}$, where $\Sa_{\bar i}$ is the
subspace of
elements of parity $\bar i$.
\begin{proposition}\label{P:knsuper}\cite{Schlsa}
The space $\Sa$ with the above introduced parity and product
is a Lie superalgebra.
\end{proposition}
\begin{definition}
The algebra $\Sa$ is the 
Krichever-Novikov type Lie superalgebra (associated to $(\Sigma,A)$).
\end{definition}
Classically this Lie superalgebra corresponds 
 to the  Neveu-Schwarz superalgebra. 
See in this context also
 \cite{Bryant},  \cite{BMRR}, \cite{BMTW}.

\subsection{Jordan Superalgebra}

Leidwanger and Morier-Genoux introduced in \cite{LeiMor}  a
Jordan superalgebra in the Krichever-Novikov setting, i.e.
\begin{equation}
\mathcal{J}:=\Fl[0]\oplus\Fl[-1/2]=
\mathcal{J}_{\bar 0}\oplus
\mathcal{J}_{\bar 1}.
\end{equation}
Recall that $\A=\Fl[0]$ is the associative algebra of 
meromorphic functions.
They define the (Jordan) product $\circ$ is  via the algebra structures 
 for the spaces $\Fl$ by
\begin{equation}
\begin{aligned}[]
f\circ g&:=f\;\sbul\; g\quad \in\Fl[0],
\\
f\circ\varphi&:=f\;\sbul\; \varphi \quad \in\Fl[-1/2]
\\
\varphi\circ\psi&:=[\varphi,\psi]\quad \in\Fl[0].
\end{aligned}
\end{equation}
By rescaling the second definition with the factor 1/2 one obtains
a Lie antialgebra. See \cite{LeiMor} for more details
and additional results on representations.
\subsection{Current Algebras}

We start with 
$\g$  a complex finite-dimensional Lie algebra and  endow
the tensor product $\gb=\g\otimes_\C \A$ with the Lie bracket
\begin{equation}
[x\otimes f, y\otimes g]=[x,y]\otimes f\cdot g,
\qquad  x,y\in\g,\quad f,g\in\A.
\end{equation} 
The algebra  $\gb$ is the higher genus current algebra.
It is an infinite dimensional Lie algebra and might be considered 
as the Lie algebra of $\g$-valued meromorphic functions on the
Riemann surface with poles only outside of $A$.
Note that we allow also the case of $\g$ an abelian Lie algebra. 
\begin{definition}
The algebra $\gb$ is called \emph{current algebra} of Krichever
Novikov type (associated to $(\Sigma,A)$).
\end{definition}
Sometimes also the name \emph{loop algebra} is used.

In the classical case 
 the current algebra $\gb$ is
the standard current algebra
$\gb=\g\otimes \C[z^{-1},z]$  with
Lie bracket
\begin{equation}
[x\otimes z^n, y\otimes z^m]=[x,y]\otimes z^{n+m}
\qquad  x,y\in\g,\quad n,m\in\Z.
\end{equation}
To point out the dependence on the geometrical structure we 
added always ``(associated to $(\Sigma,A)$) in the definition.
For simplicity we will drop it starting from now.

\section{Almost-Graded Structure}
\label{sec:almgrad}
\subsection{Definition of Almost-Gradedness}
Recall the  classical situation. This is   
the Riemann surface  $\Pro^1(\C)=S^2$, i.e. the Riemann surface 
of genus zero,
 and the points where
poles are allowed are $\{0,\infty\}$). In this case
 the algebras introduced in the last
section are graded algebras. 
In the higher genus case and even in the genus zero case with more
than
two points where poles are allowed there is no non-trivial grading
anymore.
As realized by Krichever and Novikov \cite{KNFa}
there is a weaker concept, an almost-grading, which
to  a large extend is a valuable replacement of a honest grading.
Such an almost-grading is induced by a splitting of the set $A$ into
two non-empty and disjoint sets $I$ and $O$. 
The almost-grading is fixed by exhibiting certain basis elements 
in the spaces $\Fl$ as homogeneous.
\begin{definition}\label{D:almgrad} 
Let $\La$ be a Lie or an associative algebra such that
$\La=\oplus _{n\in\Z}\La_n$ is a vector space direct sum, then
$\La$ is called an \emph{almost-graded} 
(Lie-) algebra if
\begin{enumerate}
\item[(i)] $\dim \La_n<\infty$,
\item[(ii)]
There exists constants $L_1,L_2\in\Z$ such that 
\begin{equation}
\La_n\cdot \La_m\subseteq \bigoplus_{h=n+m-L_1}^{n+m+L_2}
\La_h,\qquad \forall n,m\in\Z.
\end{equation}
\end{enumerate}
The elements in $\La_n$ are called {\it homogeneous} elements of
degree $n$, and $\La_n$ is called \emph{homogeneous subspace}
 of degree $n$.

If $\dim\La_n$ is bounded with a bound independent of $n$ we call $\La$ 
\emph{strongly almost-graded}. If we drop the condition that
$\dim\La_n$ is finite we call $\La$ \emph{weakly almost-graded}.
\end{definition}

In a similar manner almost-graded modules over almost-graded algebras
are defined. We can extend 
in an obvious way the definition to
superalgebras, resp. even to more general algebraic structures.
This definition makes complete sense also for more general index sets
$\J$. In fact we will consider the index set
$\J=(1/2)\Z$ in the case of superalgebras.
The even elements (with respect to the super-grading) will have 
integer degree, the odd elements half-integer degree.

\subsection{Separating Cycle and Krichever-Novikov Duality}
\label{sec:kndual}
Let $C_i$ be  positively oriented (deformed) circles around
the points $P_i$ in $I$, $i=1,\ldots, K$ 
 and $C_j^*$ positively oriented ones around
the points $Q_j$ in $O$, $j=1,\ldots, M$.

A cycle $C_S$ is called a separating cycle 
if it is smooth, positively oriented of multiplicity one and if 
it separates
the in- from the out-points. It might have multiple components. 
In the following we will integrate  
meromorphic differentials on $\Sigma$ without poles in 
$\Sigma\setminus A$ over closed curves $C$. 
Hence, we might consider  the $C$ and $C'$ as
equivalent if $[C]=[C']$ in  $\Ho(\Sigma\setminus A,\Z)$.
In this sense
we can write 
for every separating cycle
\begin{equation}\label{eq:cs}
[C_S]=\sum_{i=1}^K[C_i]=-\sum_{j=1}^M [C^*_j].
\end{equation}
The minus sign appears due to the opposite orientation.
Another way for giving such a $C_S$ is  via
level lines of a ``proper time
evolution'', for which I refer to Ref.~\cite{Schlce}.

Given such a separating cycle $C_S$ (resp. cycle class) 
we define a linear map
\begin{equation}
\Fl[1]\to\C,\qquad \omega\mapsto \cins \omega.
\end{equation}
As explained above the map will not depend on 
the separating line $C_S$ chosen, as
two of such will be homologous and the poles of $\omega$ are only located in
$I$ and $O$.

Consequently,
the integration of $\omega$ over $C_S$ can 
also be described
over the 
special cycles $C_i$ or equivalently over $C_j^*$.
This integration 
corresponds to calculating residues
\begin{equation}\label{eq:res}
\omega\quad\mapsto\quad
\cins \omega\ =\ \sum_{i=1}^K\res_{P_i}(\omega)
\ =\ -\sum_{l=1}^{M}\res_{Q_l}(\omega).
\end{equation}
\begin{definition}
The pairing 
\begin{equation}\label{eq:kndual}
\Fl\times\Fl[1-\la]\to \C,\quad
(f,g)\mapsto\kndual{f}{g}:=\cins f\cdot g,
\end{equation}
between $\la$ and $1-\la$ forms is called
{\it Krichever-Novikov (KN) pairing}.
\end{definition}
Note that the pairing depends not only  on $A$ (as the $\Fl$ depend on
it) but also critically on the splitting of $A$ into $I$ and $O$
as the integration path will depend on it. Once the splitting
is fixed the pairing will be fixed too.

By exhibiting dual basis elements further down we will see that
it is non-degenerate.

\subsection{The Homogeneous Subspaces}
Depending on whether $\la$ is integer or half-integer we set
$\Jl=\Z$ or $\Jl=\Z+1/2$.
For $\Fl$ we introduce for $m\in\Jl$ 
subspaces $\Fl_m$ of dimension $K$, where
$K=\#I$, 
by exhibiting  certain elements  $f_{m,p}^\la\in \Fl$,
$ p=1,\ldots, K$ which constitute a  basis of   $\Fl_m$.
Recall that the spaces $\Fl$ for $\lambda\in\Z+1/2$ depend
on the chosen square root $L$ (the 
theta characteristic) of $\K$.
The elements are the elements of degree $m$.
As explained in the following, the 
degree is in an essential way related to 
the zero orders of the elements at the points in $I$.

Let $I=\{P_1,P_2,\ldots, P_K\}$ 
then we have for the zero-order at the point $P_i\in I$
of the element  $f^\la_{n,p}$
\begin{equation}
\ord_{P_i}(f^\la_{n,p})
=(n+1-\la)-\delta_{i}^p,\quad i=1,\ldots, K\; .
\end{equation}
The prescription at the points in $O$ is made in such a way
that the element  $f^\la_{m,p}$ is essentially uniquely given.
Essentially unique means up to multiplication with a constant%
\footnote{Strictly speaking, there are some special cases where
some constants  have to be added such  that the Krichever-Novikov
duality \refE{knd} is valid, see \cite{Schlce}.}.
After fixing as additional geometric data a system of coordinates
$z_l$ centered at $P_l$ for $l=1,\ldots, K$ 
and requiring that 
\begin{equation}\label{eq:fnorm}
f_{n,p}^\la(z_p)=z_p^{n-\la}(1+O(z_p))(dz_p)^\la
\end{equation}
the element $f_{n,p}$ is uniquely fixed.
In fact, the  element $f_{n,p}^\la$ only depends on the first
jet of the coordinate $z_p$ \cite{SSpt}.

\begin{example}
Here we will not give the general recipe for the
prescription at the points in $O$, see 
\cite{Schlce}, \cite{SchlDiss}, \cite{Schlknbook}.
Just to give an example which is also an
important special case, assume $O=\{Q\}$ is a one-element
set. If either the genus $g=0$, or  $g\ge 2$, $\la\ne 0,\, 1/2,\, 1$ 
 and the points in
$A$ are in generic position 
then we require     
\begin{equation}
\ord_{Q}(f^\la_{n,p})
=-K\cdot(n+1-\la)+(2\la-1)(g-1).
\end{equation}
In the other cases (e.g. for $g=1$) there are 
some modifications at the point in $O$ necessary
for finitely
many $n$.
\end{example}

\begin{theorem}\label{T:basis}
\cite{Schlce}, \cite{SchlDiss}, \cite{Schlknbook}
Set
\begin{equation}\label{eq:kbasis}
\mathcal{B}^\la:=\{\, 
 f_{n,p}^\la\mid n\in\Jl,\  p=1,\ldots, K\,\}.
\end{equation}
Then 
(a) $\mathcal{B}^\la$ is a basis of the vector space $\Fl$.

\noindent
(b) 
The introduced basis $\mathcal{B}^\la$ of  
$\Fl[\la]$ and  $\mathcal{B}^{1-\la}$ 
of $\Fl[1-\la]$ 
are dual to each other with respect to the
Krichever-Novikov pairing \refE{kndual}, i.e.  
\begin{equation}\label{eq:knd}
\kndual{f_{n,p}^\la}{f_{-m,r}^{1-\la}}=
\delta_p^r\;\delta_n^m, \quad \forall n,m\in\Jl,\quad 
r,p=1,\ldots, K.
\end{equation}
\end{theorem}
{}From part (b) of the theorem it follows that the  
Krichever-Novikov pairing is non-degenerate. Moreover,
any element $v\in\Fl[1-\la]$ acts as linear form on $\Fl$ via
\begin{equation}
 \Phi_v:  \Fl\mapsto \C,\quad
w\mapsto \Phi_v(w):=\kndual{v}{w}.
\end{equation}
Via this pairing $\Fl[1-\la]$ can be considered as subspace of
$({\Fl})^*$. But I like to stress the fact that the identification depends
on the splitting of $A$ into $I$ and $O$ as the
KN pairing depends on it.

The full space $({\Fl})^*$ can even be described with the help
of the pairing. Consider the series
\begin{equation}\label{eq:infin}
\hat v:=\sum_{m\in\Z}\sum_{p=1}^K a_{m,p}f_{m,p}^{1-\la} 
\end{equation}
as a formal series, then
$\Phi_{\hat v}$ (as a distribution) is a well-defined  element of
${\Fl}^*$, as it will be only evaluated for finitely many
basis elements in $\Fl$.
Vice versa, every element of  ${\Fl}^*$ can be given by a suitable
$\hat v$. 
Every $\phi\in(\Fl)^*$ is uniquely given by the scalars 
$\phi(f_{m,r}^\la)$.
We set 
\begin{equation}
\hat v:=\sum_{m\in\Z}\sum_{p=1}^K \phi(f_{-m,p}^\la)\,f_{m,p}^{1-\la}. 
\end{equation}
Obviously, $\Phi_{\hat v }=\phi$.
For more information about this ``distribution
interpretation''
see \cite{SchlDiss}, \cite{SchlHab}.

The dual elements of $\La$ 
will be given by the formal series 
\refE{infin}  with basis elements from $\Fl[2]$ 
the quadratic differentials, 
the dual elements of $\A$ correspondingly from $\Fl[1]$
the differentials, 
and the
dual elements of $\Sp$ correspondingly from   $\Fl[3/2]$.
The spaces $\Fl[2]$, $\Fl[1]$ and $\Fl[3/2]$ themselves can be
considered as some kind of restricted duals.

It is quite convenient to use special notations for elements
of some important weights:
\begin{equation}
\begin{gathered}
e_{n,p}:=f_{n,p}^{-1},
\quad
\varphi_{n,p}:=f_{n,p}^{-1/2},
\quad
A_{n,p}:=f_{n,p}^{0},
\\
\omega^{n,p}:=f_{-n,p}^{1},\quad
\Omega^{n,p}:=f_{-n,p}^{2}.
\end{gathered}
\end{equation}
In view of \refE{knd} for the forms of weight 1 and 2 it is
convenient to invert the index $n$ and write it as a superscript.

\subsection{The Algebras}
\begin{theorem}
\cite{Schlce}, \cite{SchlDiss}, \cite{Schlknbook} 
\label{T:coeff}
There exists constants $R_1$ and  $R_2$ 
(depending on the genus $g$, and on  the number and splitting of the points
in $A$)  
independent of 
$n,m\in\J$ such that for the basis elements  
\begin{equation}\label{eq:coeff}
\begin{aligned}[]
f_{n,p}^\la\,\sbul\, f_{m,r}^\nu=&
\quad f_{n+m,r}^{\la+\nu}\delta_p^r
\\ &\qquad\qquad
+\sum_{h=n+m+1}^{n+m+R_1}\sum_{s=1}^Ka_{(n,p)(m,r)}^{(h,s)}f_{h,s}^{\la+\nu},
\qquad a_{(n,p)(m,r)}^{(h,s)}\in\C,
\\ \hbox{}
\\
[f_{n,p}^\la, f_{m,r}^\nu]=&
\qquad
(-\la m+\nu n)\,f_{n+m,r}^{\la+\nu+1}\delta_p^r
\\
&\qquad\qquad +
\sum_{h=n+m+1}^{n+m+R_2}\sum_{s=1}^Kb_{(n,p)(m,r)}^{(h,s)}f_{h,s}^{\la+\nu+1},
\quad b_{(n,p)(m,r)}^{(h,s)}\in\C.
\end{aligned}
\end{equation}
\end{theorem}
This says in particular that with respect to both the associative and 
Lie structure the algebra $\mathcal{F}$ is weakly almost-graded.
In  generic situations and for $N=2$ points one obtains
$R_1=g$ and $R_2=3g$.

The reason why we only have weakly almost-gradedness is that 
\begin{equation}
\Fl=\bigoplus_{m\in\Jl}\Fl_m,\qquad\text{with}\quad \dim \Fl_m=K.
\end{equation}
If we add up for a fixed $m$ all $\la$ we get that our
homogeneous spaces are infinite dimensional.

In the definition of our KN type algebra 
only finitely many $\la$ are involved, hence the following is
immediate
\begin{theorem}\label{T:almgrad}
The Krichever-Novikov type 
vector field algebras $\La$, function algebras $\A$, differential 
operator algebras $\Do$, Lie superalgebras $\Sa$, and 
Jordan superalgebras $\Ja$ are all (strongly) almost-graded.
\end{theorem} 
We obtain
\begin{equation}
\dim\La_n=\dim\A_n=K,\quad
\dim\Sa_n=\dim\Ja_n=2K,\quad
\dim\Do_n=3K.
\end{equation}
If $\mathcal{U}$ is one  of these algebras, with product
denoted by $[\;,\;]$ then
\begin{equation}
[\mathcal{U}_n,\mathcal{U}_m]
\subseteq \bigoplus_{h=n+m}^{n+m+R_i}\mathcal{U}_h,
\end{equation}
with $R_i=R_1$ for $\mathcal{U}=\A$ and $R_i=R_2$ otherwise.

For further reference let us specialize the lowest degree
term component in \refE{coeff} for certain special cases.
\begin{equation}
\begin{aligned}[]
A_{n,p}\cdot A_{m,r}\ &=\ A_{n+m,r}\,\delta_r^p\ +\ \text{h.d.t.}
\\
A_{n,p}\cdot f_{m,r}^\la\ &=\ f_{n+m,r}^\la\,\delta_r^p\ +\ \text{h.d.t.}
\\
[e_{n,p},e_{m,r}]\ &=\ (m-n)\cdot e_{n+m,r}\,\delta_r^p\ +\ \text{h.d.t.}
\\
e_{n,p}\ldot f_{m,r}^\la\ &= \ (m+n\la)\cdot 
f_{n+m,r}^\la\,\delta_r^p\ +\ \text{h.d.t.}
\end{aligned}
\end{equation}
Here $\ \text{h.d.t.}\ $ denote linear combinations of 
basis elements of degree between
\newline
$n+m+1$ and $n+m+R_i$,

Finally, the almost-grading of $\A$ induces an almost-grading of
the current algebra $\ \gb\ $ by setting
$\gb_n=\g\otimes \A_n$.
We obtain
\begin{equation}
\gb=\bigoplus_{n\in\Z}\gb_n,\quad \dim\gb_n=K\cdot\dim\g.
\end{equation}
\subsection{Triangular Decomposition and Filtrations}
Let $\U$ be one of the above introduced algebras (including the
current algebra). 
On the basis of the almost-grading we obtain a triangular
decomposition of the algebras
\begin{equation} 
\U=\U_{[+]}\oplus\U_{[0]}\oplus\U_{[-]},
\end{equation}
where 
\begin{equation}
\U_{[+]}:=\bigoplus_{m>0}\U_m,\quad
\U_{[0]}=\bigoplus_{m=-R_i}^{m=0}\U_m,\quad
\U_{[-]}:=\bigoplus_{m<-R_i}\U_m.
\end{equation}
By the almost-gradedness the $[+]$ and
$[-]$ subspaces are  
(infinite dimensional) subalgebras. The $[0]$ spaces in general
not. Sometimes we will use {\it critical strip} for $\U_{[0]}$.

With respect to the almost-grading 
of $\Fl$ we can introduce a filtration
\begin{equation}
\begin{gathered}
\Fl_{(n)}:=\bigoplus_{m\ge n} \Fl_m,
\\ ....\quad
\supseteq\quad  \Fl_{(n-1)}\quad
\supseteq \quad\Fl_{(n)}\quad
\supseteq \quad \Fl_{(n+1)}\quad ....
\end{gathered}
\end{equation}
\begin{proposition}
\begin{equation}
\Fl_{(n)}:=\{\ f\in\Fl\mid \ord_{P_i}(f)\ge n-\la,\forall i=1,\ldots, K\
\}.
\end{equation}
\end{proposition}
This proposition is very important. 
In case that $O$ has more than one point
there are certain choices, e.g. numbering of the points in
$O$, different rules, etc. involved in defining the almost-grading. 
Hence, if the choices are made differently the subspaces $\Fl_n$ 
might depend on them, and consequently also the almost-grading.
But by this proposition the induced filtration is indeed
canonically defined  via the splitting of $A$ into $I$ and $O$.

Moreover, different choices will give equivalent almost-grading.
We stress the fact, that under a KN algebra we will always
understand one of introduced the algebras (or even some
others still to come) together with an almost-grading (resp.
equivalence class of almost-grading) introduced by the
splitting $A=I\cup O$.

\section{Central Extensions}
\label{sec:centvec}

Central extension of our algebras appear naturally in the context of
quantization and regularization of actions. In \refS{seminf} we will
see a typical example. Of course they are also of independent 
mathematical interest.

\subsection{Central Extensions and Cocycles}

A central extension of a Lie algebra $W$ is a special
Lie algebra structure on the 
vector space direct sum $\widehat{W}=\C\oplus W$.
If we denote  $\hat x:=(0,x)$ and $t:=(1,0)$ then the Lie structure is
given by 
\begin{equation}\label{eq:cext}
[\hat x, \hat y]=\widehat{[x,y]}+\psi(x,y)\cdot t,
\quad [t,\widehat{W}]=0,
\quad x,y\in W.
\end{equation}
The map  $x\mapsto \hat x=(0,x)$ is a linear splitting map.
$\widehat{W}$ will be a Lie algebra, e.g. will fulfill the 
Jacobi identity, if and only if $\psi$ is antisymmetric and fulfills
the Lie algebra 2-cocycle condition
\begin{equation}
0=d_2\psi(x,y,z):=
\psi([x,y],z)+
\psi([y,z],x)+
\psi([z,x],y).
\end{equation}
A 2-cochain $\psi$ is a coboundary
if  there exists 
a $\varphi:W\to\C$ such that 
\begin{equation}
\psi(x,y)=\varphi([x,y]).
\end{equation}
One easily shows that a coboundary is a cocycle. Hence, we can define
the second Lie algebra cohomology $\Ho^2(W,\C)$ 
of $W$ with 
values in the trivial module $\C$ as this quotient. 

There is the notion of equivalence of central extensions. 
For the definition in terms of short exact sequences, I refer to the
standard text books. 
Equivalently, they can be described by a different choice of the
linear splitting map. Instead of  
\newline
$x\mapsto \hat x=(0,x)$ one
chooses  $x\mapsto {\hat x}'=(\varphi,x)$ with $\varphi:W\to\C$ a
linear form.
For $\varphi\ne 0$ the Lie structure corresponding to 
\refE{cext} will be different, but equivalent.
In fact in \refE{cext} we obtain a different 2-cocycle
$\psi'$.

{}From this definition follows  that 
that two central extensions are equivalent if and only if the difference of 
their defining
2-cocycles $\psi$ and $\psi'$ is a coboundary.
In this way the second Lie algebra cohomology $\Ho^2(W,\C)$ 
classifies equivalence classes
of central extensions. The class $[0]$ corresponds to the
trivial central extension. In this case the splitting map
is a Lie homomorphism. 
To construct central extensions of our algebras we have to 
find such Lie algebra 2-cocycles. 

Clearly, equivalent central extensions are isomorphic.
The opposite is not true. 
Furthermore, in our case we can always rescale
the central element by multiplying it with a nonzero scalar. 
This is not an equivalence of
central extensions but nevertheless an irrelevant modification.
Hence we will mainly be interested in central extensions modulo 
equivalence and rescaling. They are classified by 
$[0]$ and the elements of the projectivized cohomology 
space $\Pro( \Ho^2(W,\C))$.

In the classical case we have $\dim\Ho^2(\W,\C)=1$, hence there are
only two essentially different central extensions, the 
splitting one given by the direct sum $\C\oplus \W$ of Lie
algebras and the up to equivalence and rescaling unique non-trivial
one, the Virasoro algebra $\V$.
\begin{remark}
Given a vector space bases $\{e_\rho\mid \rho\in\mathcal{R}\}$ of
$W$, a vector space basis of $\widehat{W}$ will be given 
by $\{\hat e_\rho:=(0,e_\rho)\mid \rho\in\mathcal{R}, t:=(1,0)\;\}$.
An equivalent central extension can be described as a change
of basis and rescaling of the form
\begin{equation}
{\hat e_\rho}'=\hat e_\rho+\varphi( e_\rho)\,t,
\qquad
t'=\alpha\cdot t,\ \alpha\in\C^*.
\end{equation}
\end{remark}

\subsection{Geometric Cocycles}\label{sec:geocyc}
The defining cocycle for the Virasoro algebra obviously does not
make any sense in the higher genus and/or multi-point case. We need 
a geometric description. For this we have first to introduce
connections.

\subsubsection{Projective and Affine Connections}

Let $\ (U_\al,z_\al)_{\al\in J}\ $ be a covering of the Riemann
surface
$\Sigma$ 
by holomorphic coordinates with transition functions
$z_\be=f_{\be\al}(z_\al)$.

\begin{definition}
(a) A system of local (holomorphic, meromorphic) functions 
$\ R=(R_\al(z_\al))\ $ 
is called a (holomorphic, meromorphic) {\it projective 
connection} if it transforms as
\begin{equation}\label{eq:pc}
R_\be(z_\be)\cdot (f_{\beta,\alpha}')^2=R_\al(z_\al)+S(f_{\beta,\alpha}),
\qquad\text{with}\quad
S(h)=\frac {h'''}{h'}-\frac 32\left(\frac {h''}{h'}\right)^2,
\end{equation}
the Schwartzian derivative.
Here ${}'$ denotes differentiation with respect to
the coordinate $z_\al$.

(b) A system of local (holomorphic, meromorphic) functions 
$\ T=(T_\al(z_\al))\ $ 
is called a (holomorphic, meromorphic) {\it affine 
connection} if it transforms as
\begin{equation}\label{eq:zc}
T_\be(z_\be)\cdot (f_{\beta,\alpha}')=T_\al(z_\al)+
\frac{f''_{\beta,\alpha}}{f'_{\beta,\alpha}}.
\end{equation}
\end{definition}

Every Riemann surface admits a holomorphic projective 
connection  \cite{HawSchiff}, \cite{Gun}.
Given a point $P$ then there exists always a meromorphic
 affine connection holomorphic outside of $P$ and having maximally a
 pole
of order one there \cite{SchlDiss}.

 From their very definition 
  it follows that the difference of two affine (projective) connections
 will be a (quadratic) differential.
 Hence, after fixing one affine (projective) connection all others are obtained
 by adding (quadratic) differentials.

\subsubsection{The Function Algebra $\A$}
We consider $\A$ as an abelian Lie algebra. Let $C$ be an arbitrary smooth
but not necessarily connected curve. We set
\begin{equation}\label{eq:exta}
\psi^1_{C}(g,h):=\cint{C} gdh,\quad g,h\in\A.
\end{equation}
\subsubsection{The Current Algebra $\gb$}
For $\gb=\g\otimes \A$ we first have to 
fix  $\beta$ a symmetric, invariant, bilinear form on 
$\g$ (not necessarily 
non-degenerate).
Invariance means  that we have 
$\beta([x,y],z)=\beta(x,[y,z])$ for all $x,y,z\in\g$.
The cocycle is given as
\begin{equation}\label{eq:daff}
\psi^2_{C,\beta}(x\otimes g, y\otimes h):=
\beta(x,y)\cdot \cint{C} gdh,\quad  x,y\in\g,\ g,h\in\A.
\end{equation}
\subsubsection{The Vector Field Algebra  $\La$}
Here it is a little bit more delicate. First we have to choose 
a (holomorphic) projective connection $R$.
We define 
\begin{equation}\label{eq:vecg}
\psi^3_{C,R}(e,f):=\cintl{C} \left(\frac 12(e'''f-ef''')
-R\cdot(e'f-ef')\right)dz\ .
\end{equation}
Only by the term related with the projective connection it will be
a well-defined differential, i.e. independent of the coordinate
chosen.
It is shown in \cite{SchlDiss} that it 
is a cocycle. Another choice of a projective connection
will result in a cohomologous one.
Hence, the equivalence class of the central extension will be the same.

\subsubsection{The Differential Operator Algebra  $\Do$}
For the differential operator algebra the cocycles for $\A$ can be
extended by zero on the subspace $\La$. 
The cocycles for $\La$ can be pulled back.
In addition there is a third type of cocycles mixing $\A$ and $\La$:
\begin{equation}\label{eq:mix}
\psi^4_{C,T}(e,g):=
\cintl{C}(e g''+T eg')dz,
\qquad  e\in\La,g\in\A,
\end{equation} 
with an affine connection $T$, with at most a pole of order one
at a fixed point in $O$.
Again, a different choice of the connection will not change
the equivalence class.
For more details on the cocycles see \cite{Schlloc}.

\subsubsection{The Lie Superalgebra $\Sa$} 
Here we have to take into account that it is not a Lie algebra. Hence,
the Jacobi identity has to be replaced by the super-Jacobi identity.
The conditions for being a cocycle for the superalgebra
cohomology will change too.
Recall the definition of the algebra from \refS{super},
in particular that the even elements (parity $0$) are the vector
fields
and the odd elements (parity $1$) are the half-forms.
A bilinear form $c$ is a cocycle if the following is true.
The bilinear map $c$  will be symmetric if $x$ and $y$ are odd, otherwise
it will be antisymmetric.
\begin{equation}\label{eq:csymm}
c(x,y)=-(-1)^{\bar x \bar y}c(x,y).
\end{equation}
The super-cocycle condition reads as 
\begin{equation}\label{eq:scocyc}
(-1)^{\bar x \bar z}c(x,[y,z])+
(-1)^{\bar y \bar x}c(y,[z,x])+
(-1)^{\bar z \bar y}c(z,[x,y])\ = \ 0.
\end{equation}
With the help of $c$ we can define central extensions in the Lie
superalgebra
sense. If we put the condition that the central element is even
then the cocycle $c$ has to be an even map and $c$ vanishes 
for pairs of elements of different parity. 

By  convention we  denote vector fields by $e,f,g,...$ and
-1/2-forms by $\varphi,\psi,\chi,..$ and get
\begin{equation}
c(e,\varphi)=0,\quad e\in\La,\ \varphi\in\Sp.
\end{equation}
The super-cocycle conditions for the even elements is just the
cocycle condition for the Lie subalgebra $\La$. The only other
nonvanishing super-cocycle condition is for the 
{\it (even,odd,odd)} elements and reads as
\begin{equation}\label{eq:cyccon}
c(e,[\varphi,\psi])-c(\varphi,e\ldot \psi)-c(\psi,e\ldot \varphi)
=0.
\end{equation}
Here the definition of the product $[e,\psi]:=e\ldot\psi$
was used.

If we have a cocycle $c$ for the algebra $\Sa$ we obtain
by restriction a cocycle for the algebra $\La$.
For the mixing term we know that $c(e,\psi)=0$.
A naive try to put just anything for $c(\varphi,\psi)$ 
(for example 0) will not work
as \refE{cyccon} relates the restriction of the cocycle on $\La$ with
its values on $\Sp$.
\begin{proposition}\cite{Schlsa}
Let $C$ be any closed (differentiable) curve on $\Sigma$ not meeting
the
points in $A$, and let $R$ be any (holomorphic) projective connection
then the bilinear extension of 
\begin{equation}
\label{eq:sacoc}
\begin{aligned}
\Phi_{C,R}(e,f)&:=\cintl{C} \left(\frac 12(e'''f-ef''')
-R\cdot(e'f-ef')\right)dz
\\
\Phi_{C,R}(\varphi,\psi)
&:=
-\cintl{C} \left(\varphi''\cdot \psi+\varphi\cdot \psi''
-R\cdot\varphi\cdot \psi\right)dz
\\
\Phi_{C,R}(e,\varphi)&:=0
\end{aligned}
\end{equation}
gives a Lie superalgebra cocycle for $\Sa$, hence defines a central
extension of $\Sa$.
A different projective connection will yield a cohomologous 
cocycle.
\end{proposition}
A similar formula was given  by
Bryant in \cite{Bryant}. By adding the projective connection in
the second part of 
\refE{sacoc} he corrected some formula appearing in 
\cite{BMRR}. He only considered the two-point case and only the
integration over a separating cycle.
See also \cite{Kreusch} for the multi-point case, where  still only
the integration over a separating cycle is considered.

In contrast to the differential operator algebra case the two parts 
cannot be prescribed independently. Only with the same integration path
(more precisely, homology class) and the given factors in front
of the integral it will work.
The reason for this that  \refE{cyccon} relates both.

\subsection{Uniqueness and Classification of Central Extensions}
\label{sec:class}

Our cocycles depend on the choice of the connections 
$R$ and $T$.
But different choices will not change the equivalence class. Hence,
this ambiguity does not disturb us.
What really matters is that they depend on the integration curve $C$ 
chosen.

In contrast to the classical situation, 
for the higher genus and/or multi-point situation
there are many essentially
different closed curves and hence many non-equivalent central
extensions defined by the integration.

But we should take into account that we want to extend the
almost-grading from our algebras to the centrally extended ones.
This means we take $\deg\hat x:=\deg x$ and assign a degree 
$deg(t)$ to the central element $t$, and still obtain
almost-gradedness.

This is possible if and only if our defining cocycle $\psi$ is 
``local'' in the following sense (the name was introduced 
in the two point case by
Krichever and Novikov in \cite{KNFa}). There exists 
$M_1,M_2\in\Z$ such that
\begin{equation}
\forall n,m:\quad\psi(W_n,W_m)\ne 0\ \implies\  
M_1\le n+m\le M_2.
\end{equation}
Here $W$ stands for any  of our algebras (including the supercase).
Very important, ``local'' is defined in terms of the
almost-grading, and the grading itself depends on the splitting
$A=I\cup O$. Hence what is ``local''  depends on the splitting too.

We will call a cocycle {\it bounded} (from above) if there exists
$M\in\Z$ such that
\begin{equation}
\forall n,m:\quad\psi(W_n,W_m)\ne 0\ \implies\  
n+m\le M.
\end{equation}
Similarly bounded from below is defined. Locality means bounded
from above and below.

Given a cocycle class we call it bounded (resp. local) if and only if
it contains a representing cocycle which is bounded (resp. local).
Not all cocycles in a bounded class have to be bounded.
If we choose as integration path a separating cocycle $C_S$,
or one of the $C_i$  then
the above introduced geometric cocycles  are local,
resp. bounded.
Recall that in this case integration can  be done by calculating
residues at the in-points or  at the out-points.
All these cocycles are cohomologically nontrivial.
The theorems in the following  concern the opposite direction.
They were treated in my works \cite{Schlloc},
\cite{Schlaff}, \cite{Schlsa}.

I  start with the vector field case as it will give a model for
all the results.

\begin{theorem}\label{T:vclass}
\cite{Schlloc} 
Let $\La$ be the Krichever--Novikov vector field algebra.
\newline
(a) The space of bounded 
cohomology classes is $K$-dimensional ($K=\#I$).
A basis is given by setting  the integration path in
\refE{vecg}  to $C_i$, $i=1,\ldots, K$ the little (deformed) circles
around the points $P_i\in I$.
\newline
(b) The space of local cohomology classes is one-dimensional.
A generator is given by integrating \refE{vecg} over a
separating cocycle $C_S$.
\newline
(c) Up to equivalence and rescaling there is only one 
non-trivial one-dimensional
central extension of the vector field algebra $\La$ which 
allows  an extension of the almost-grading.
\end{theorem}
Part (c) says that the result  is completely analogous to the case 
of the Witt algebra.
Here I like to repeat again  the fact that for $\La$ depending on the set
$A$ and its possible splittings into two disjoint subsets
there are different almost-gradings. Hence, the ``unique'' central
extension finally obtained will also depend on the splitting.
Only in the two point case there is only one splitting
possible.

The above theorem is a model for all other classification
results. We will always obtain a statement about the bounded
(from above) cocycles and then for the local cocycles.

As $\A$ is an abelian Lie algebra every skew-symmetric bilinear form
will be a non-trivial cocycle. Hence, there is no hope of uniqueness.
But if we add the condition of
$\La$-invariance, i.e.
\begin{equation}
\psi(e.g,h)+\psi(g,e.h)=0,\quad
\forall e\in\La,\  g,h\in\A
\end{equation}
things will change. 

Let us denote the the subspace of local cohomology classes
by $\Ho^2_{loc}$, and the subspace of local and $\La$-invariant by
$\Ho^2_{\La,loc}$. Note that the condition is only required for at
least one representative in the cohomology class.
We collect a part of the results for the other algebras 
in the following
theorem.
\begin{theorem}
$ $

\smallskip
\noindent
(1)\quad  $\dim\Ho^2_{\La,loc}(\A,\C)=1$,

\smallskip
\noindent
(2)\quad  $\dim\Ho^2_{loc}(\Sa,\C)=1$,

\smallskip
\noindent
(3)\quad  $\dim\Ho^2_{loc}(\Do,\C)=3$,

\smallskip
\noindent
(4)\quad  $\dim\Ho^2_{loc}(\gb,\C)=1$\quad  for 
$\g$ a simple finite-dimensional
Lie algebra,

\smallskip
\noindent

\noindent
A basis of the cohomology spaces are given by taking the cohomology
classes
of the cocycles 
\refE{exta}, \refE{daff}, \refE{vecg}, \refE{mix}, \refE{sacoc}
obtained by integration over a separating cycle $C_S$.
\end{theorem}
Correspondingly, we obtain also for these algebras the corresponding
result about uniqueness of almost-graded central extensions.
For the differential operator algebra we got three independent
cocycles. This generalizes results of \cite{ACKP} for the
classical case.

For the bounded cocycle classes we have
to multiply the dimensions above by $K$.
For 
the supercase with odd central elements 
the bounded cohomology vanishes. 

For $\g$ a reductive Lie algebra and the cocycle $\La$-invariant
if restricted to the abelian part,
a complete classification  of local cocycle classes 
for $\gb$ can be found 
in \cite{Schlaff}.
Note that in the case of a simple Lie algebra every 
symmetric, invariant bilinear form $\beta$ is a multiple 
of the Cartan-Killing form.

I like to mention that in all the applications I know of, the cocycles
coming
from representations, regularizations, etc. are local.
Hence, the uniqueness or classification can be used.

\section{Further Results}\label{sec:fr}

Above the basic concepts, results about the structure
of these Krichever-Novikov type algebras and their central extensions
were treated.
Of course, this does not close the story. 
I will add some further important constructions and applications 
but due to space limitations
only in a very
condensed 
manner.

\subsection{Semi-Infinite Forms and Fermionic 
Fock Space Representations}\label{sec:seminf}
$ $

Our Krichever-Novikov  vector field algebras $\La$ have as Lie modules
the spaces $\Fl$. These representations are not of the type 
physicists
are usually interested in, as  there is no ground state (no vacuum).
There are neither  creation nor annihilation  operators 
which can be used to
construct the full representation out of a vacuum state. 

To obtain such desired representations the almost-grading comes into
play.
First, using the grading of $\Fl$  it is possible to construct 
starting from  $\Fl$, the forms of weight $\la\in 1/2\Z$,
the
{\it semi-infinite wedge forms} $\Hl$.

The vector space $\Hl$
 is generated by basis elements which are formal expressions
\begin{equation}\label{eq:sform}
\Phi=f^{\la}_{(i_1)}\wedge f^{\la}_{(i_2)} \wedge 
f^{\la}_{(i_3)}\wedge \cdots,
\end{equation}
where $(i_1)=(m_1,p_1)$ is a double index indexing our basis elements.
The indices are in strictly increasing  lexicographical
order. They are stabilizing in the sense that they will increase
exactly by one starting from a certain point, depending on $\Phi$.
The action of $\La$ can be extended by Leibniz rule from $\Fl$ to
$\Hl$.
But a problem arises. For elements of the critical strip 
$\La_{[0]}$ of the
algebra $\La$ it might happen that it produces infinitely many
contributions.
The action has to be regularized (as physicists like to call it,
but it is a well-defined mathematical procedure). 

Here the almost-grading has his second appearance.
By the (strong) almost-graded module structure 
of $\Fl$ the algebra $\La$ can be
imbedded into the Lie algebra of both-sided infinite matrices
\begin{equation}
\glib:=\{A=(a_{ij})_{i,j\in\Z}\mid
\exists r=r(A), \text{ such that }
a_{ij}=0 \text{ if } |i-j|>r\ \},
\end{equation}
with ``infinitely many diagonals''.
The embedding will depend on the weight $\la$.
For $\glib$ there exists a procedure for the regularization of
the action on the semi-infinite wedge forms
\cite{DJKM}, \cite{KaP}, see also \cite{KaRa} for 
a nice pedagogical treatment. 
In particular, there is a unique non-trivial central extension
$\glih$. If we pull-back the defining cocycle for the 
extension 
we obtain a central extension $\Lh_{\la}$ of $\La$ and 
the required regularization of the action  of   $\Lh_{\la}$ on
$\Hl$.
As the embedding of $\La$ depends on the weight $\la$ the 
cocycle will do so.
The pull-back cocycle will be local.
 Hence, 
by the classification results of \refS{class} it is the
unique central extension class defined by \refE{vecg}
integrated over $C_S$ (up to a $\la$ dependent rescaling).

In $\Hl$ there are invariant subspaces, which are
generated by a certain ``vacuum vectors''. 
Such a vacuum $\Phi_T$ is given by an element of the 
form \refE{sform} which starts with the element $f^{\la}_{(T,1)}$,
and the indices for  the following ones increase always by one.
The subalgebra 
$\La_{[+]}$ annihilates the vacuum, the central element and the 
other elements of degree zero act by multiplication with a
constant on the vacuum and the whole representation state is generated by
$\La_{[-]}\oplus \La_{[0]}$ from the vacuum.

As the function algebra $\A$ operates as multiplication operators
on $\Fl$ the above representation can be extended to the algebra 
$\Do$ (see details in
\cite{SchlDiss}) after one passes over to central extensions.
The cocycle again is local and hence, up to coboundary,
 it will be 
a certain linear combination of the 3 generating cocycles for the
differential operator algebra.
In fact it will be 
\begin{equation}\label{eq:dhl}
c_\la[\psi^3_{C_S}]+\frac {2\la-1}{2}[\psi^4_{C_S}]-[\psi^1_{C_S}],
\qquad c_\la:=-2(6\la^2-6\la+1).
\end{equation}
Recall that  $\psi^3$ is the cocycle for the vector field algebra,
$\psi^1$ the cocycle for the function algebra,
and $\psi^4$ the mixing cocycle, see \cite{Schlloc} for
details.
Note that the expression for $c_\la$ appears also in Mumford's 
formula \cite{SchlRS}
relating divisors on the moduli space of curves.

Also the representation on $\Hl$ gives a 
projective representation of the algebra
of $\Dal$ of differential operators of all orders. 
It is exactly the combination \refE{dhl} which lifts to a cocycle for 
 $\Dal$ and gives a central extension 
 $\Dalh$.
For $\La$ we could rescale the central element. Hence the 
central extension $\Lh$ did not depend essentially on 
the weight. Here this is different.
The central extension $\Dh_{\la}$ depends on it.

\smallskip

For the centrally extended current algebras $\gh$,
the affine algebra of KN type,
in a similar way fermionic Fock space representations can be
constructed, see \cite{Shf}, \cite{SSpt}.

\subsection{$b$ -- $c$ Systems}

Related to the above there are other quantum algebra systems which
can be realized on $\Hl$.
On the space $\Hl$ the forms $\Fl$ act by wedging  elements
$f^\la\in\Fl$ in front of the semi-infinite wedge form,
i.e.
\begin{equation}
\Phi\mapsto  f^\la\wedge \Phi.
\end{equation}
Using the Krichever-Novikov duality
pairing \refE{kndual} to contract in the
semi-infinite wedge form the entries 
with the  form $f^{1-\la}\in\Fl[1-\la]$ the latter form will act
$\Hl$.
For  $\Phi$  a basis element \refE{sform} of $\Hl$
the contraction is defines via
\begin{equation}
i(f^{1-\la})\Phi
=
\sum_{l=1}^\infty
(-1)^{l-1}\kndual{f^{1-\la}}{f^\la_{i_l}}
\cdot f^{\la}_{(i_1)}\wedge f^{\la}_{(i_2)} \wedge \cdots
\check {f}^{\la}_{(i_l)}\cdots.
\end{equation}
Here $\check {f}^{\la}_{(i_l)}$ indicates as usual that this  element
will not be there anymore.

Both operations together create a Clifford algebra structure, which is
sometimes called a $\ b-c\ $ system, see
\cite{SchlDiss}, \cite{SchlHab}, \cite{Schlknbook}.

\subsection{Sugawara Representation}

Given  an admissible  representation
of the 
centrally extended current algebra $\gh$ we can construct  the so called 
{\it Sugawara operators}. Here admissible means, that the
central element operates as constant $\times$ identity, and that
every element $v$ in the representation space will be annihilated by
the  elements in $\gh$ of sufficiently high degree (the degree 
depends on the
element $v$).  
The Sugawara operator is an infinite formal sum of operators and
is constructed as  the product of the current
operators which are again formal infinite sum of operators.
To make the product well-defined a normal ordering has to be set,
which moves the annihilation operators to the right to act first.
It turns out that after some rescaling the operators  appearing
in the formal sum of the Sugawara operators 
give a representation of a centrally extended vector field algebra
$\La$.
The central extension is due to 
the appearance of the normal ordering.
Again the defining cocycle is local and we know that the central
extension defined by the representation is the central extension
given by our geometric cocycle $\psi^3_{C_S}$.
See \cite{SchlShSug}, \cite{SchlSug}, \cite{SchlHab}
for details.
\subsection{Wess-Zumino-Novikov-Witten models and
  Knizhnik-Zamolodchikov  Connection}

Despite the fact, that it is a very important application,
the following description is extremely condensed.
More 
can be found in \cite{SSpt}, \cite{SSpt2}, 
\cite{Schlglob}.
See also  \cite{SheBo}, \cite{Schlknbook}.
Wess-Zumino-Novikov-Witten (WZNW)  models are defined on the basis of 
a fixed finite-dimensional Lie algebra $\g$. One considers families of
representations of the affine KN algebras $\gh$ 
(which is an almost-graded central extension of the current algebra
$\gb$
of KN type) 
defined over the
moduli space of Riemann surfaces of genus $g$ with $K+1$ marked points
and splitting of type $(K,1)$.
The single point in $O$ will be a reference point.
The data of the moduli of the Riemann surface and the marked points
enter the definition of the algebra $\gh$ and the representation.
The construction of certain co-invariants
 yields  a special vector bundle of finite rank over 
moduli space, called the vector bundle of conformal blocks.
With the help of
the Krichever Novikov vector field algebra, and using the
Sugawara construction,
 the {\it Knizhnik-Zamolodchikov
(KZ)  connection} is given. It is projectively  flat.
An essential fact is that certain elements in the critical strip
$\La_{[0]}$ of the vector field algebra correspond to infinitesimal 
deformations of the moduli and to moving the marked points.
This gives a global operator approach in contrast to the semi-local
approach of Tsuchiya, Ueno, and Yamada \cite{TUY}.

\subsection{Geometric Deformations of the Witt and 
Virasoro Algebra}\label{sec:geodef}

As the second Lie algebra cohomology of the Witt and Virasoro algebra 
with values in
their adjoint module vanishes \cite{FiaSchl1},
\cite{Schlrig}, \cite{Fiarig}
both are formally and infinitesimally rigid. This means that all
formal (and infinitesimal) families 
where the  special fiber is one of 
these algebras are equivalent to the trivial 
one.
Nevertheless, we showed in \cite{FiaSchl1}
that there exists naturally defined families of Krichever-Novikov
vector field algebras 
defined for the torus with two marked points \cite{SchlDeg}, 
\cite{Deck}, \cite{RDS}.
These families are obtained by a geometric degeneration process.
The families 
have as special element the Witt algebra (resp. Virasoro
algebra).
All other fibers are non-isomorphic to it.
Hence, these families are even locally non-trivial.
This is a phenomena which can only be observed for infinite
dimensional algebras. See also the case of affine Lie algebra and
some general treatment in
 \cite{FiaSchlaff}, \cite{FiaSchlaffo}, \cite{Schlwvc}.

\section{Lax Operator Algebras}
\label{sec:lax}

Recently, a new class of current type algebras appeared, the
Lax operator algebras.
As the naming indicates,
they are related to infinite dimensional integrable systems
\cite{ShLax}.
The algebras  were introduced by Krichever \cite{Klax}, and Krichever and
Sheinman \cite{KSlax}.
Here I will report on their definition.

Compared to the KN current type algebra we will allow 
additional singularities which will play a special role.
The points where these  singularities are 
allowed are called {\it weak singular points}. The set is denoted
by 
\begin{equation}
W=\{\gamma_s\in\Sigma\setminus A\mid s=1,\ldots, R\}.
\end{equation}
Let $\g$ be  one of the classical matrix
algebras $\gl(n)$, $\sln(n)$, $\so(n)$, $\spn(2n)$.
We assign to every point $\ga_s$ 
 a vector $\al_s\in\C^n$ 
(resp.  $\in \C^{2n}$ for  $\spn(2n)$).
The system 
\begin{equation}
\mathcal{T}:=\{(\gamma_s,\al_s)\in\Sigma\times \C^n\mid s=1,\ldots, R\}
\end{equation}
is called 
{\it Tyurin data}. 
\begin{remark}
In case that $R=n\cdot g$ 
and for  generic values of $(\ga_s,\al_s)$ with $\al_s\ne 0$ the tuples
of pairs
 $(\ga_s,[\al_s])$ with $[\al_s]\in\Pro^{n-1}(\C)$
parameterize semi-stable rank $n$ and degree $n\, g$ framed holomorphic
vector bundles as shown by Tyurin \cite{Tyvb}.
Hence, the name Tyurin data. 
\end{remark}
We consider $\g$-valued meromorphic functions
\begin{equation}
L:\ \Sigma\ \to\  \g,
\end{equation}
which are holomorphic outside  $W\cup A$, have at most
poles of order one (resp. of order two for $\spn(2n)$) at the
points in $W$, and fulfill certain conditions at $W$ depending on
$\mathcal{T}$. 
To describe them let us fix local coordinates
$w_s$ centered at $\ga_s$,  $s=1,\ldots, R$.
For {\bf $\gl(n)$} the conditions 
are as follows.
For $s=1,\ldots, R$ we require that there exist $\be_s\in\C^n$ 
and $\ka_s\in \C$ such that the
function $L$ has the following expansion at $\ga_s\in W$ 
\begin{equation}\label{eq:glexp}
L(w_s)=\frac {L_{s,-1}}{w_s}+
L_{s,0}+\sum_{k>0}L_{s,k}{w_s^k},
\end{equation}
with
\begin{equation}\label{eq:gldef}
L_{s,-1}=\al_s \be_s^{t},\quad
\tr(L_{s,-1})=\be_s^t \al_s=0,
\quad
L_{s,0}\,\al_s=\ka_s\al_s.
\end{equation}
In particular, if $L_{s,-1}$ is 
non-vanishing then it is a rank 1 matrix, and if 
$\al_s\ne 0$  then it is  
an eigenvector of $L_{s,0}$.
The requirements \refE{gldef} are independent of the chosen
coordinates $w_s$.

It is not at all clear that the commutator of two such 
matrix functions fulfills again these conditions. 
But it is shown in \cite{KSlax} that they indeed close to
a Lie algebra (in fact in the case of $\gl(n)$ they constitute
an associative algebra under the matrix product).
If one of the $\alpha_s=0$ then 
the conditions at the point $\gamma_s$ correspond to the fact, 
that $L$ has to be
holomorphic there. 
If all $\alpha_s$'s are zero or $W=\emptyset$ we obtain back
the current algebra of KN type. In some sense the Lax operator
algebras generalize them. But their definition is  restricted to the
case that our finite-dimensional Lie algebra has to be one 
from the list.
In the bundle interpretation of the Tyurin data the
KN case corresponds to the trivial rank $n$ bundle.

For {\bf $\sln(n)$} the only 
additional condition  is that 
in \refE{glexp} all matrices $L_{s,k}$ are  trace-less.
The conditions \refE{gldef}  remain unchanged.

In the case of {\bf $\so(n)$} one requires that
all $L_{s,k}$ in \refE{glexp} are  skew-symmetric.
In particular, they are trace-less.
Following \cite{KSlax} the set-up has to be slightly modified.
First only  those Tyurin parameters $\al_s$ are allowed which satisfy
$\al_s^t\al_s=0$.
Then the 1. relation in \refE{gldef} is changed to obtain
\begin{equation}\label{eq:sodef}
L_{s,-1}=\al_s\be_s^t-\be_s\al_s^t,
\quad
\tr(L_{s,-1})=\be_s^t\al_s=0,
\quad
L_{s,0}\,\al_s=\ka_s\al_s.
\end{equation}

\medskip
For {\bf $\spn(2n)$}
we consider  a symplectic form  $\hat\sigma$  
for $\C^{2n}$ given by
a non-degenerate skew-symmetric matrix $\sigma$.
The Lie algebra $\spn(2n)$ is the Lie algebra of
matrices $X$ such that $X^t\sigma+\sigma X=0$.
The condition $\tr(X)=0$ will be automatic. 
At the weak singularities we have the expansion
\begin{equation}\label{eq:glexpsp}
L(w_s)=\frac {L_{s,-2}}{w_s^2}+\frac {L_{s,-1}}{w_s}+
L_{s,0}+L_{s,1}{w_s}+\sum_{k>1}L_{s,k}{w_s^k}.
\end{equation}
The condition \refE{gldef} is  modified as
follows (see \cite{KSlax}):
there exist $\be_s\in\C^{2n}$, 
$\nu_s,\ka_s\in\C$ such that 
\begin{equation}\label{eq:spdef}
L_{s,-2}=\nu_s \al_s\al_s^t\sigma,\quad
L_{s,-1}=(\al_s\be_s^t+\be_s\al_s^t)\sigma,
\quad{\be_s}^t\sigma\al_s=0,\quad
L_{s,0}\,\al_s=\kappa_s\al_s.
\end{equation}
{}Moreover, we require
$\al_s^t\sigma L_{s,1}\al_s=0.
$
Again under the point-wise
matrix commutator the set of such maps constitute a Lie algebra.

\bigskip
The next step is  to introduce an almost-graded structure for these
Lax operator algebras induced by a splitting of the set
$A=I\cup O$.  
This is done for the two-point case in \cite{KSlax} and for
the multi-point case in \cite{SSlaxm}. 
{}From the applications there is again a need to classify 
almost-graded central
extensions.  The author obtained this jointly with O. Sheinman
in \cite{SSlax}, \cite{SSlaxm},
see also \cite{Schllax} for an overview.
For the Lax operator algebras associated to 
the simple algebras $\sln(n),\so(n), \spn(n)$ 
it will be unique (meaning: given a splitting of $A$ there is
an almost-grading and with respect to this there is up to 
equivalence and rescaling only  one non-trivial almost-graded
central extension).
For $\gl(n)$ we obtain two independent local cocycle classes if
we assume $\La$-invariance on the reductive part.

%
%
%

\section{Some Historical Remarks}
\label{sec:hist}

In this section I will give some historical remarks 
(also on related works)
on 
Krichever-Novikov type algebras. Space limitations do not allow to
give a complete list of references. For this I have to refer to
\cite{Schlknbook}.

In 1987 the ground-breaking work of Krichever and Novikov
\cite{KNFa}, \cite{KNFb}, \cite{KNFc} in the two-point case
initiated the subject. 
They introduced the vector field algebra, the function algebra and
the affine algebra with their almost-graded structure.
To acknowledge their work these algebras are nowadays
called Krichever-Novikov (KN) type algebras. Sheinman joint in by
investigating the affine algebras and their representations
\cite{Shea},
\cite{Shd},
\cite{Sha},
\cite{Shdd},
\cite{Shf},
\cite{Shcas}.

As it should have become clear from this review what should be called
a KN type algebra is not the algebra alone but the algebra together
with a chosen almost-grading.
It is exactly the step of introducing such an almost-grading
which is not straightforward.

{}From the application in CFT, string theory, etc. there was clearly the
need to pass over to a multi-point picture.
In the multi-point case this was done 1989 by the author in 
\cite{Schlmp}, \cite{Schleg}, \cite{Schlce}, \cite{SchlDiss}.
The almost-grading is  induced by splitting of the set  $A$
of points where poles are allowed 
into two non-empty disjoint subsets $I$ and $O$.
In the applications quite often such a splitting is naturally given.
In the context of fields and strings the points correspond to
incoming and outgoing fields and strings respectively.
Without considering an almost-grading Dick \cite{Di}
gave also a generalization of the vector field algebra. He obtained 
results similar to \cite{Schlmp}.
Only in the work of Sadov \cite{Sad} 1990 an almost-grading 
is also discussed.

Note that in the two-point case there is only one
splitting. Hence, quite often one does not mention explicitely the
grading for the Witt and Virasoro algebra, respectively 
the almost-grading for the KN type algebras.
Nevertheless, the grading is heavily used.

The genus zero and two-point case is the classical well-studied case.
But already for genus zero and more than two points there are 
interesting things so study, see 
\cite{Schleg},
\cite{FiaSchl1}, \cite{FiaSchlaff},
\cite{Brem2}, \cite{Brem3}, \cite{Schlmunich}.
For genus one, the complex torus case, there is
\cite{Schleg}, \cite{Brem1}, \cite{Deck}, \cite{RDS}, 
\cite{FiaSchl1}, \cite{FiaSchlaff}.

After the work of Krichever and Novikov appeared physicists  got very 
much interested in these algebras and the possibilities of using
these objects for a global operator approach. Especially I like
to mention the work of the people around Bonora
\cite{BMRR},
\cite{BLMR},
\cite{BMTW},
\cite{BRRW} and by 
\cite{Bryant}.
This includes also  the superversions.
A lot of more names could be given.

It is not only the possible applications in physics which
makes the KN type algebras  so interesting. From the mathematical point 
general infinite dimensional Lie algebras are hard to approach.
KN type algebras supply examples of them which are given in a
geometrical
context, hence (hopefully) better to understand.
A typical example of this are the families of geometric deformation
of the Witt algebra which I mentioned in \refS{geodef}
obtained by degenerations of tori.
Quite recently also in the context of Jordan Superalgebras and
Lie antialgebras \cite{Ovs1} examples    
were constructed on the basis of KN type algebras
\cite{LeiMor}, \cite{LeiMor1}, \cite{Kreusch}

In this review we 
discussed extensively the case of 
2nd Lie algebra cohomology with values in the trivial module.
But we did not touch the question of the 
general Lie algebra 
cohomology of Krichever-Novikov type algebras.
Here I refer e.g. to work of Wagemann 
\cite{Wag1}, \cite{Wag2}.


\end{document}